\documentclass[reqno]{amsart}
\usepackage{amssymb,hyperref}

\theoremstyle{plain}
\newtheorem{thm}{Theorem}
\newtheorem{lem}{Lemma}
\newtheorem{cor}{Corollary}

\theoremstyle{definition}
\newtheorem{ex}{Example}

\renewcommand{\Re}{\mathrm{Re}}

\title
[Properties of certain analytic functions associated with ...]
{Properties of certain analytic functions \\
associated with two boundary points}

\author{Hitoshi Shiraishi}
\address{Hitoshi Shiraishi \newline
Department of Mathematics \newline
Kinki University \newline
Higashi-Osaka, Osaka 577-8502, Japan}
\email{shiraishi@math.kindai.ac.jp}

\author{Shigeyoshi Owa}
\address{Shigeyoshi Owa \newline
Department of Mathematics \newline
Kinki University \newline
Higashi-Osaka, Osaka 577-8502, Japan}
\email{owa@math.kindai.ac.jp}


\date{}

\begin{document}

\begin{abstract}
For analytic functions $f(z)$ in the closed unit disk $\overline{\mathbb{U}}$,
two boundary points $z_1$ and $z_2$ such that $\alpha = (f'(z_1)+f'(z_2))/2 \in f'(\mathbb{U})$ are considered.
The object of the present paper is to discuss some interesting conditions for $f(z)$ to be $|f'(z)-1|<\rho|1-\alpha|$ in $\mathbb{U}$ with some examples.
\end{abstract}

\begin{flushleft}
This paper was published in the journal: \\
Int. J. Math. Anal. (Ruse) {\bf 4} (2010), No. 25-28, 1271--1284. \\
{\footnotesize \url{http://www.m-hikari.com/ijma/ijma-2010/ijma-25-28-2010/owaIJMA25-28-2010.pdf}}
\end{flushleft}
\hrule

\

\

\maketitle

\section{Introduction}

\

Let $\mathcal{A}_{n}$ denote the class of functions
$$
f(z)=z+a_{n+1}z^{n+1}+a_{n+2}z^{n+2}+ \ldots
\qquad(n=1,2,3,\ldots)
$$
that are analytic in the closed unit disk $\overline{\mathbb{U}}=\{z \in \mathbb{C}:|z|\leqq1\}$
and $\mathcal{A}=\mathcal{A}_1$.
Also,
the open unit disk is denoted by $\mathbb{U}=\{z \in \mathbb{C}:|z|<1\}$.

For two boundary points $z_1$ and $z_2$, let us consider
$$
\alpha = \frac{f'(z_1)+f'(z_2)}{2} \in f'(\mathbb{U}).
$$

For such $\alpha$ $(\alpha \ne 1)$, if $f(z) \in \mathcal{A}_n$ satisfies
$$
\left|\frac{f'(z)-\alpha}{1-\alpha} -1\right| < \rho
\qquad(z \in \mathbb{U})
$$
for some real $\rho > 0$, then
$$
|f'(z)-1| < \rho|1-\alpha|
\qquad(z \in \mathbb{U}).
$$

Therefore, if $0 < \rho|1-\alpha| < 1$, then $f(z)$ is close-to-convex (univalent) in $\mathbb{U}$. In the present paper, we use such a technique for $f(z)$.

\

The basic tool in proving our results is the following lemma due to Jack \cite{m1ref1}
(also, due to Miller and Mocanu \cite{m1ref2}).

\

\begin{lem} \label{jack} \quad
Let the function $w(z)$ defined by
$$
w(z)=a_nz^n+a_{n+1}z^{n+1}+a_{n+2}z^{n+2}+ \ldots
\qquad(n=1,2,3,\ldots)
$$
be analytic in $\mathbb{U}$ with $w(0)=0$.
If $\left|w(z)\right|$ attains its maximum value on the circle $\left|z\right|=r$ at a point $z_{0}\in\mathbb{U}$,
then there exists a real number $k \geqq n$ such that
$$
\frac{z_{0}w'(z_{0})}{w(z_{0})}=k
$$
and
$$
\Re\left(\frac{z_0w''(z_0)}{w'(z_0)}\right)+1\geqq k.
$$
\end{lem}

\

\section{Main results}

\

Applying Lemma \ref{jack},
we drive the following results.

\

\begin{thm} \label{d4thm1} \quad
If $f(z)\in\mathcal{A}_{n}$ satisfies
$$
\left| \frac{zf''(z)}{f'(z)} \right|
< \frac{|1-\alpha|n\rho}{1+|1-\alpha|\rho}
\qquad(z\in\mathbb{U})
$$
for some complex $\alpha=\dfrac{f'(z_1)+f'(z_2)}{2}\in f'(\mathbb{U})$ and $\alpha\neq1$
such that $z_1,z_2 \in \partial \mathbb{U}$,
and for some real $\rho>1$,
then
$$
|f'(z)-1|<\rho|1-\alpha|
\qquad(z\in\mathbb{U}).
$$
\end{thm}

\

\begin{proof}\quad
Let us define $w(z)$ by
\begin{align}
w(z)
&= \frac{f'(z)-\alpha}{1-\alpha}-1
\qquad(z\in\mathbb{U})
\label{d4thm1eq1}\\
&= \frac{(n+1)a_{n+1}}{1-\alpha}z^n+\frac{(n+2)a_{n+2}}{1-\alpha}z^{n+1}+ \ldots . \nonumber
\end{align}

Then, clearly, $w(z)$ is analytic in $\mathbb{U}$ and $w(0)=0$.
Differentiating both sides in (\ref{d4thm1eq1}),
we obtain
$$
\frac{zf''(z)}{f'(z)}
= \frac{(1-\alpha)zw'(z)}{(1-\alpha)w(z)+1},
$$
and therefore,
$$
\left| \frac{zf''(z)}{f'(z)} \right|
= \left| \frac{(1-\alpha)zw'(z)}{(1-\alpha)w(z)+1} \right|
< \frac{|1-\alpha|n\rho}{1+|1-\alpha|\rho}
\qquad(z\in\mathbb{U}).
$$

If there exists a point $z_0 \in \mathbb{U}$ such that
$$
\max_{|z| \leqq |z_{0}|} |w(z)|
= |w(z_{0})|
= \rho,
$$
then Lemma \ref{jack} gives us that $w(z_0)=\rho e^{i\theta}$ and $z_0 w'(z_0)=kw(z_0)$ ($k \geqq n$).

For such a point $z_0$,
we have
\begin{align*}
&\left| \frac{z_0f''(z_0)}{f'(z_0)} \right| \\
&= \left| \frac{(1-\alpha)z_0w'(z_0)}{(1-\alpha)w(z_0)+1} \right| \\
&= \left| \frac{(1-\alpha)kw(z_0)}{(1-\alpha)w(z_0)+1} \right| \\
&= \frac{|(1-\alpha)|k\rho}{|(1-\alpha)\rho e^{i\theta}+1|} \\
&= \frac{|(1-\alpha)|k\rho}{\sqrt{(1+|1-\alpha|\rho\cos\phi)^2+(|1-\alpha|\rho\sin\phi)^2}} \\
&= \frac{|(1-\alpha)|k\rho}{\sqrt{1+|1-\alpha|^2\rho^2+2|1-\alpha|\rho\cos\phi}} \\
&\geqq \frac{|1-\alpha|k\rho}{1+|1-\alpha|\rho} \\
&\geqq \frac{|1-\alpha|n\rho}{1+|1-\alpha|\rho}.
\end{align*}
where $\phi=\theta+\arg(1-\alpha)$.

This contradicts our condition in the theorem.
Therefore, there is no $z_0 \in \mathbb{U}$ such that $|w(z_0)|=\rho$.
This means that $|w(z)|<\rho$ for all $z \in \mathbb{U}$.
It follows that
$$
\left| \frac{f'(z)-\alpha}{1-\alpha}-1 \right|
< \rho
\qquad(z\in\mathbb{U})
$$
so that
$\left| f'(z)-1 \right|<\rho|1-\alpha|$ in $\mathbb{U}$.
\end{proof}

\

\begin{ex} \label{d4ex1} \quad
Let us consider a function
$$
f(z)
= z+a_{n+1}z^{n+1}
\qquad(z\in\mathbb{U})
$$
with $|a_{n+1}|<\dfrac{1}{2(n+1)}$.
Differenciating the function $f(z)$,
we obtain
$$
\frac{zf''(z)}{f'(z)}
=\frac{n(n+1)a_{n+1}z^n}{1+(n+1)a_{n+1}z^n},
$$
and therefore,
\begin{align*}
\left| \frac{zf''(z)}{f'(z)} \right|
&= \left| \frac{n(n+1)a_{n+1}z^n}{1+(n+1)a_{n+1}z^n} \right| \\
&< \frac{n(n+1)|a_{n+1}|}{1-(n+1)|a_{n+1}|} \qquad (z \in \mathbb{U}).
\end{align*}

We choose two boundary points $z_1$ and $z_2$ such that
$f'(z_1)=1+(n+1)|a_{n+1}|$ and $f'(z_2)=1+(n+1)|a_{n+1}|i$.
Then we know that
$ z_1=e^{-i\frac{\arg(a_{n+1})}{n}}$ and $ z_2=e^{i\frac{\pi-2\arg(a_{n+1})}{2n}} $.
For such $z_1$ and $z_2$, we obtain
$$
\alpha
= \frac{f'(z_1) + f'(z_2)}{2} = 1 + \frac{(n+1)|a_{n+1}|(1+i)}{2}
$$
so that
$$
1 - \alpha
= -\,\frac{(n+1)|a_{n+1}|(1+i)}{2}.
$$
Now, we consider some $\rho > 1$ such that
$$
\frac{n(n+1)|a_{n+1}|}{1-(n+1)|a_{n+1}|}
\leqq \frac{|1-\alpha|n\rho}{1+|1-\alpha|\rho}
= \frac{n(n+1)|a_{n+1}|\rho}{\sqrt{2}+(n+1)|a_{n+1}|\rho}.
$$

This gives us that
$$
\rho
\geqq \frac{\sqrt{2}}{1-2(n+1)|a_{n+1}|}.
$$

For such $\alpha$ and $\rho$, we know that $f(z)$ satisfies
$$
|f'(z) -1|
< (n+1)|a_{n+1}|
\leqq \frac{(n+1)|a_{n+1}|}{1-2(n+1)|a_{n+1}|}
\leqq \rho|1-\alpha|.
$$
\end{ex}

\

If we defined
\begin{equation}
\max_{z \in \mathbb{U}} |f'(z)-\alpha|
= M_\alpha,
\label{d4cor1eq1}
\end{equation}
then we can have
\begin{align*}
|f'(z)-1|
&< \rho |1-\alpha| \\
&= \rho |f'(0)-\alpha| \\
&\leqq \rho M_\alpha
\qquad(z\in\mathbb{U}).
\end{align*}

So, we get

\

\begin{cor} \label{d4cor1.1} \quad
If $f(z)\in\mathcal{A}_{n}$ satisfies
$$
\left| \frac{zf''(z)}{f'(z)} \right|
< \frac{|1-\alpha|n\rho}{1+|1-\alpha|\rho}
\qquad(z\in\mathbb{U})
$$
for some complex $\alpha=\dfrac{f'(z_1)+f'(z_2)}{2}\in f'(\mathbb{U})$ and $\alpha\neq1$
such that $z_1,z_2 \in \partial \mathbb{U}$,
and for some real $\rho>1$,
then
$$
|f'(z)-1|<\rho M_\alpha
\qquad(z\in\mathbb{U}).
$$
\end{cor}

\

Putting
\begin{equation}
\alpha
= \frac{f'(z_1)+f'(z_2)+ \ldots +f'(z_m)}{m}
\label{d4cor1eq2}
\end{equation}
for $z_1,z_2, \ldots ,z_m \in \partial \mathbb{U}$,
we obtain

\

\begin{cor} \label{d4cor1.2} \quad
If $f(z)\in\mathcal{A}_{n}$ satisfies
$$
\left| \frac{zf''(z)}{f'(z)} \right|
< \frac{|1-\alpha|n\rho}{1+|1-\alpha|\rho}
\qquad(z\in\mathbb{U})
$$
for some complex $\alpha=\dfrac{f'(z_1)+f'(z_2)+ \ldots +f'(z_m)}{m} \in f'(\mathbb{U})$
such that \\
$z_1,z_2, \ldots ,z_m \in \partial \mathbb{U}$ and $\alpha\neq1$,
and for some real $\rho>1$,
then
$$
|f'(z)-1|<\rho|1-\alpha|
\qquad(z\in\mathbb{U}).
$$
\end{cor}

\

We also derive

\

\begin{thm} \label{d4thm2} \quad
If $f(z)\in\mathcal{A}_{n}$ satisfies
$$
\left| zf''(z)-\frac{zf''(z)}{f'(z)} \right|
< \frac{|1-\alpha|^2 n \rho^2}{1+|1-\alpha|\rho}
\qquad(z\in\mathbb{U})
$$
for some complex $\alpha=\dfrac{f'(z_1)+f'(z_2)}{2}\in f'(\mathbb{U})$ and $\alpha\neq1$
such that $z_1,z_2 \in \partial \mathbb{U}$,
and for some real $\rho>1$,
then
$$
|f'(z)-1|<\rho|1-\alpha|
\qquad(z\in\mathbb{U}).
$$
\end{thm}

\

\begin{proof}\quad
Define $w(z)$ by
\begin{align}
w(z)
&= \frac{f'(z)-\alpha}{1-\alpha}-1
\qquad(z\in\mathbb{U})
\label{d4thm2eq1}\\
&= \frac{(n+1)a_{n+1}}{1-\alpha}z^n+\frac{(n+2)a_{n+2}}{1-\alpha}z^{n+1}+ \ldots. \nonumber
\end{align}

Evidently, $w(z)$ is analytic in $\mathbb{U}$ and $w(0)=0$.
Differentiating (\ref{d4thm2eq1}) logarithmically and simplyfing,
we have
\begin{align*}
zf''(z)-\frac{zf''(z)}{f'(z)}
&=zf''(z) \left( 1-\frac{1}{f'(z)} \right) \\
&= \frac{(1-\alpha)^2 zw'(z)w(z)}{(1-\alpha)w(z)+1},
\end{align*}
and hence,
$$
\left| zf''(z)-\frac{zf''(z)}{f'(z)} \right|
= \left| \frac{(1-\alpha)^2 zw'(z)w(z)}{(1-\alpha)w(z)+1} \right|
< \frac{|1-\alpha|^2 n \rho^2}{1+|1-\alpha|\rho}
\qquad(z\in\mathbb{U}).
$$

If there exists a point $z_0 \in \mathbb{U}$ such that
$$
\max_{|z| \leqq |z_{0}|} |w(z)|
= |w(z_{0})|
= \rho,
$$
then Lemma \ref{jack} gives us that $w(z_0)=\rho e^{i\theta}$ and $z_0 w'(z_0)=kw(z_0)$ ($k \geqq n$).

For such a point $z_0$,
we have
\begin{align*}
&\left| z_0f''(z_0)-\frac{z_0f''(z_0)}{f'(z_0)} \right| \\
&= \left| \frac{(1-\alpha)^2 z_0w'(z_0)w(z_0)}{(1-\alpha)w(z_0)+1} \right| \\
&= \frac{|1-\alpha|^2 \rho^2 k}{|(1-\alpha)\rho e^{i\theta}+1|} \\
&= \frac{|1-\alpha|^2 \rho^2 k}{\sqrt{(1+|1-\alpha|\rho\cos\phi)^2+(|1-\alpha|\rho\sin\phi)^2}} \\
&= \frac{|1-\alpha|^2 \rho^2 k}{\sqrt{1+|1-\alpha|^2\rho^2+2|1-\alpha|\rho\cos\phi}} \\
&\geqq \frac{|1-\alpha|^2 k \rho^2}{1+|1-\alpha|\rho} \\
&\geqq \frac{|1-\alpha|^2 n \rho^2}{1+|1-\alpha|\rho}.
\end{align*}
where $\phi=\theta+\arg(1-\alpha)$.

This contradicts our condition in the theorem.
Therefore, there is no $z_0 \in \mathbb{U}$ such that $|w(z_0)|=\rho$.
This means that $|w(z)|<\rho$ for all $z \in \mathbb{U}$.
This implies that
$$
\left| \frac{f'(z)-\alpha}{1-\alpha}-1 \right|
< \rho
\qquad(z\in\mathbb{U}).
$$
\end{proof}

\

\begin{ex} \label{d4ex2} \quad
We consider a function $f(z)$ given by
$$
f(z)
=z+a_{n+1}z^{n+1}
\qquad(z\in\mathbb{U})
$$
with $|a_{n+1}|< \dfrac{1}{n+1}$.
Differenciating $f(z)$,
we obtain that
$$
zf''(z)-\frac{zf''(z)}{f'(z)}
= \frac{n(n+1)^2 a_{n+1}^2 z^{2n}}{1+(n+1)a_{n+1}z^n},
$$
that is, that
\begin{align*}
\left| zf''(z)-\frac{zf''(z)}{f'(z)} \right|
&= \left| \frac{n(n+1)^2 a_{n+1}^2 z^{2n}}{1+(n+1)a_{n+1}z^n} \right| \\
&< \frac{n(n+1)^2 |a_{n+1}|^2}{1-(n+1)|a_{n+1}|} \qquad(z \in \mathbb{U}).
\end{align*}

Choosing same points $z_1$ and $z_2$ in Example \ref{d4ex1},
we see that
$$
1-\alpha
= -\,\frac{(n+1)|a_{n+1}|(1+i)}{2}.
$$

If we consider some $\rho > 1$ such that
$$
\frac{n(n+1)^2|a_{n+1}|^2}{1-(n+1)|a_{n+1}|}
\leqq \frac{|1-\alpha|^2 n\rho^2}{1+|1-\alpha|\rho},
$$
we have that
$$
\rho
\geqq \frac{\sqrt{2}}{1-(n+1)|a_{n+1}|}.
$$

For such $\alpha$ and $\rho$,
we know that
$$
|f'(z) - 1|
< (n+1)|a_{n+1}|
\leqq \frac{(n+1)|a_{n+1}|}{1-(n+1)|a_{n+1}|}
\leqq \rho|1-\alpha|.
$$
\end{ex}

\

Further, we obtain

\

\begin{cor} \label{d4cor2.1} \quad
If $f(z)\in\mathcal{A}_{n}$ satisfies
$$
\left| zf''(z)-\frac{zf''(z)}{f'(z)} \right|
< \frac{|1-\alpha|^2 n \rho^2}{1+|1-\alpha|\rho}
\qquad(z\in\mathbb{U})
$$
for some complex $\alpha=\dfrac{f'(z_1)+f'(z_2)}{2}\in f'(\mathbb{U})$ and $\alpha\neq1$
such that $z_1,z_2 \in \partial \mathbb{U}$,
and for some real $\rho> 1$,
then
$$
|f'(z)-1|<\rho M_\alpha
\qquad(z\in\mathbb{U}).
$$
\end{cor}

\

\begin{cor} \label{d4cor2.2} \quad
If $f(z)\in\mathcal{A}_{n}$ satisfies
$$
\left| zf''(z)-\frac{zf''(z)}{f'(z)} \right|
< \frac{|1-\alpha|^2 n \rho^2}{1+|1-\alpha|\rho}
\qquad(z\in\mathbb{U})
$$
for some complex $\alpha=\dfrac{f'(z_1)+f'(z_2)+ \ldots +f'(z_m)}{m} \in f'(\mathbb{U})$
such that \\
$z_1,z_2, \ldots ,z_m \in \partial \mathbb{U}$ and $\alpha\neq1$,
and for some real $\rho>1$,
then
$$
|f'(z)-1|<\rho|1-\alpha|
\qquad(z\in\mathbb{U}).
$$
\end{cor}

\

Further,
we discuss a new application for Lemma \ref{jack}.

\

\begin{thm} \label{d4thm3} \quad
If $f(z)\in\mathcal{A}_{n}$ satisfies
$$
\Re \left( \frac{z(zf''(z))'}{f'(z)-1} \right)
< n^2
\qquad(z\in\mathbb{U})
$$
for some complex $\alpha=\dfrac{f'(z_1)+f'(z_2)}{2}\in f'(\mathbb{U})$ and $\alpha\neq1$
such that $z_1,z_2 \in \partial \mathbb{U}$,
and for some real $\rho>1$,
then
$$
|f'(z)-1|<\rho|1-\alpha|
\qquad(z\in\mathbb{U}).
$$
\end{thm}

\

\begin{proof}\quad
Defining the function $w(z)$ by
\begin{align*}
w(z)
&= \frac{f'(z)-\alpha}{1-\alpha}-1
\qquad(z\in\mathbb{U}) \\
&= \frac{(n+1)a_{n+1}}{1-\alpha}z^n+\frac{(n+2)a_{n+2}}{1-\alpha}z^{n+1}+ \ldots,
\end{align*}
we have that $w(z)$ is analytic in $\mathbb{U}$ with $w(0)=0$.
Since,
$$
\frac{z(zf''(z))'}{f'(z)-1}
= \frac{zw'(z)}{w(z)}+\frac{z^2w''(z)}{w(z)},
$$
we obtain that
\begin{align*}
\Re \left( \frac{z(zf''(z))'}{f'(z)-1} \right)
&= \Re \left( \frac{zw'(z)}{w(z)}+\frac{z^2w''(z)}{w(z)} \right) \\
&= \Re \left( \frac{zw'(z)}{w(z)} \left( 1+\frac{zw''(z)}{w'(z)} \right) \right)
< n^2
\qquad(z\in\mathbb{U}).
\end{align*}

If there exists a point $z_0 \in \mathbb{U}$ such that
$$
\max_{|z| \leqq |z_{0}|} |w(z)|
= |w(z_{0})|
= \rho,
$$
then Lemma \ref{jack} gives us that $w(z_0)= \rho e^{i\theta}$ and $z_0 w'(z_0)=kw(z_0)$ ($k \geqq n$).

Thus we have
\begin{align*}
\Re \left( \frac{z_0(z_0f''(z_0))'}{f'(z_0)-1} \right)
&= \Re \left( \frac{z_0w'(z_0)}{w(z_0)} \left( 1+\frac{z_0w''(z_0)}{w'(z_0)} \right) \right) \\
&= k \left( 1+\Re \left( \frac{z_0w''(z_0)}{w'(z_0)} \right) \right) \\
&= k^2 \\
&\geqq n^2.
\end{align*}

This contradicts our condition in the theorem.
Therefore, there is no $z_0 \in \mathbb{U}$ such that $|w(z_0)|=\rho$.
This means that $|w(z)|<\rho$ for all $z \in \mathbb{U}$.
\end{proof}

\

We also have the following corollaries.

\

\begin{cor} \label{d4cor3.1} \quad
If $f(z)\in\mathcal{A}_{n}$ satisfies
$$
\Re \left( \frac{z(zf''(z))'}{f'(z)-1} \right)
< n^2
\qquad(z\in\mathbb{U})
$$
for some complex $\alpha=\dfrac{f'(z_1)+f'(z_2)}{2}\in f'(\mathbb{U})$ and $\alpha\neq1$
such that $z_1,z_2 \in \partial \mathbb{U}$,
and for some real $\rho> 1$,
then
$$
|f'(z)-1|<\rho M_\alpha
\qquad(z\in\mathbb{U}).
$$
\end{cor}

\

\begin{cor} \label{d4cor3.2} \quad
If $f(z)\in\mathcal{A}_{n}$ satisfies
$$
\Re \left( \frac{z(zf''(z))'}{f'(z)-1} \right)
< n^2
\qquad(z\in\mathbb{U})
$$
for some complex $\alpha=\dfrac{f'(z_1)+f'(z_2)+ \ldots +f'(z_m)}{m} \in f'(\mathbb{U})$
such that \\
$z_1,z_2, \ldots ,z_m \in \partial \mathbb{U}$ and $\alpha\neq1$,
and for some real $\rho>1$,
then
$$
|f'(z)-1|<\rho|1-\alpha|
\qquad(z\in\mathbb{U}).
$$
\end{cor}

\

Next our result is contained in

\

\begin{thm} \label{d4thm4} \quad
If $f(z)\in\mathcal{A}_{n}$ satisfies
$$
\frac{zf''(z)}{f'(z)-1}
\neq k
\qquad(z\in\mathbb{U})
$$
for some complex $\alpha=\dfrac{f'(z_1)+f'(z_2)}{2}\in f'(\mathbb{U})$ and $\alpha\neq1$
such that $z_1,z_2 \in \partial \mathbb{U}$,
some real $\rho>1$,
and for all real $k \geqq n$,
then
$$
|f'(z)-1|<\rho|1-\alpha|
\qquad(z\in\mathbb{U}).
$$
\end{thm}

\

\begin{proof}\quad
Let us define the function $w(z)$ by
\begin{align}
w(z)
&= \frac{f'(z)-\alpha}{1-\alpha}-1
\qquad(z\in\mathbb{U}) \\
&= \frac{(n+1)a_{n+1}}{1-\alpha}z^n+\frac{(n+2)a_{n+2}}{1-\alpha}z^{n+1}+ \ldots. \nonumber
\end{align}

Clearly, $w(z)$ is analytic in $\mathbb{U}$ with $w(0)=0$.
We want to prove that $|w(z)|<\rho$ in $\mathbb{U}$.
Note that
$$
\frac{zf''(z)}{f'(z)-1}
=\frac{zw'(z)}{w(z)}
\qquad(z\in\mathbb{U}).
$$

If there exists a point $z_0 \in \mathbb{U}$ such that
$$
\max_{|z| \leqq |z_{0}|} |w(z)|
= |w(z_{0})|
= \rho,
$$
then Lemma \ref{jack} gives us that $w(z_0)=\rho e^{i\theta}$ and $z_0 w'(z_0)=kw(z_0)$ ($k \geqq n$).

Thus we have
\begin{align*}
\frac{z_0f''(z_0)}{f'(z_0)-1}
&= \frac{z_0w'(z_0)}{w(z_0)} \\
&= k
\end{align*}

This contradicts the condition in the theorem.
Therefore, there is no $z_0 \in \mathbb{U}$ such that $|w(z_0)|=\rho$.
This means that $|w(z)|<\rho$ for all $z \in \mathbb{U}$.
We conclude that $\left| f'(z)-1 \right|<\rho|1-\alpha|$ in $\mathbb{U}$.
\end{proof}

\

\begin{cor} \label{d4cor4.1} \quad
If $f(z)\in\mathcal{A}_{n}$ satisfies
$$
\frac{zf''(z)}{f'(z)-1}
\neq k
\qquad(z\in\mathbb{U})
$$
for some complex $\alpha=\dfrac{f'(z_1)+f'(z_2)}{2}\in f'(\mathbb{U})$ and $\alpha\neq1$
such that $z_1,z_2 \in \partial \mathbb{U}$,
some real $\rho> 1$,
and for all real $k \geqq n$,
then
$$
\left| f'(z)-1 \right| < \rho M_\alpha
\qquad(z\in\mathbb{U}).
$$
\end{cor}

\

\begin{cor} \label{d4cor4.2} \quad
If $f(z)\in\mathcal{A}_{n}$ satisfies
$$
\frac{zf''(z)}{f'(z)-1}
\neq k
\qquad(z\in\mathbb{U})
$$
for some complex $\alpha=\dfrac{f'(z_1)+f'(z_2)+ \ldots +f'(z_m)}{m} \in f'(\mathbb{U})$
such that \\
$z_1,z_2, \ldots ,z_m \in \partial \mathbb{U}$,
some real $\rho>1$,
and for all real $k \geqq n$,
then
$$
|f'(z)-1|<\rho|1-\alpha|
\qquad(z\in\mathbb{U}).
$$
\end{cor}

\

Finally,
we derive

\

\begin{thm} \label{d4thm5} \quad
If $f(z)\in\mathcal{A}_{n}$ satisfies
$$
\left| \frac{zf'(z)}{f(z)}-1 \right|
< \frac{|1-\beta|n\rho}{1+|1-\beta|\rho}
\qquad(z\in\mathbb{U})
$$
for some complex $\beta=\dfrac{F(z_1) + F(z_2)}{2} \in F(\mathbb{U})$
such that $z_1,z_2 \in \partial \mathbb{U}$ and $\beta\neq1$,
and for some real $\rho>1$, where $F(z)=\dfrac{f(z)}{z}$,
then
$$
\left| \frac{f(z)}{z}-1 \right|<\rho|1-\beta|
\qquad(z\in\mathbb{U}).
$$
\end{thm}

\

\begin{proof}\quad
Let us define $w(z)$ by
\begin{align}
w(z)
&= \frac{\frac{f(z)}{z}-\beta}{1-\beta}-1
\qquad(z\in\mathbb{U})
\label{d4thm5eq1}\\
&= \frac{a_{n+1}}{1-\beta}z^n+\frac{a_{n+2}}{1-\beta}z^{n+1}+ \ldots. \nonumber
\end{align}

Then, we have that $w(z)$ is analytic in $\mathbb{U}$ and $w(0)=0$.
Differenciating (\ref{d4thm5eq1}) in both side logarithmically and simplifying,
we obtain
$$
\frac{zf'(z)}{f(z)}
= \frac{(1-\beta)zw'(z)}{(1-\beta)w(z)+1},
$$
and hence,
$$
\left| \frac{zf'(z)}{f(z)} \right|
= \left| \frac{(1-\beta)zw'(z)}{(1-\beta)w(z)+1} \right|
< \frac{|1-\beta|n\rho}{1+|1-\beta|\rho}
\qquad(z\in\mathbb{U}).
$$

Using the same process of the proof in Theorem $\ref{d4thm1}$,
we complete the proof of the theorem.
\end{proof}

\

\begin{ex} \label{d4ex5} \quad
We consider a function
$$
f(z)
=z+a_{n+1}z^{n+1}
\qquad(z\in\mathbb{U})
$$
with $|a_{n+1}| < \dfrac{1}{2}$.
Differenciating the function,
we obtain
$$
\frac{zf'(z)}{f(z)}-1
= \frac{na_{n+1}z^n}{1+a_{n+1}z^n},
$$
and therefore,
\begin{align*}
\left| \frac{zf'(z)}{f(z)}-1 \right|
&= \left| \frac{na_{n+1}z^n}{1+a_{n+1}z^n} \right| \\
&< \frac{n|a_{n+1}|}{1-|a_{n+1}|} \qquad(z \in \mathbb{U}).
\end{align*}
Consider two boundary points $z_1 = e^{-i\frac{\arg(a_{n+1})}{n}}$ and $z_2 = e^{i\frac{\pi-2 \arg(a_{n+1})}{2n}}$.
Then,
since $F(z_1)=1+|a_{n+1}|$ and $F(z_2)=1+|a_{n+1}|i$, we see that
$$
1 - \beta
= -\,\frac{|a_{n+1}|(1+i)}{2}.
$$

For such $\beta$, we consider some $\rho$ such that
$$
\frac{n|a_{n+1}|}{1-|a_{n+1}|}
\leqq \frac{|1-\beta|n\rho}{1+|1-\beta|\rho}.
$$

This gives us that
$$
\rho
\geqq \frac{\sqrt{2}}{1-2|a_{n+1}|}.
$$

Therefore,
we have that
$$
\left|\frac{f(z)}{z} - 1\right|
< |a_{n+1}|
\leqq \frac{|a_{n+1}|}{1-2|a_{n+1}|}
\leqq \rho|1-\beta|.
$$
\end{ex}

\

Defining $M_\beta$ by
$$
\max_{z \in \mathbb{U}} \left| \frac{f(z)}{z}-\beta \right|
= M_\beta,
$$
we have the following corollary.

\

\begin{cor} \label{d4cor5.1} \quad
If $f(z)\in\mathcal{A}_{n}$ satisfies
$$
\left| \frac{zf'(z)}{f(z)}-1 \right|
< \frac{|1-\beta|n\rho}{1+|1-\beta|\rho}
\qquad(z\in\mathbb{U})
$$
for some complex $\beta=\dfrac{F(z_1) + F(z_2)}{2}\in F(\mathbb{U})$
such that $z_1,z_2 \in \partial \mathbb{U}$ and $\beta\neq1$,
and for some real $\rho> 1$,
where $F(z)=\dfrac{f(z)}{z}$,
then
$$
\left| \frac{f(z)}{z}-1 \right|<\rho M_\beta
\qquad(z\in\mathbb{U}).
$$
\end{cor}

\

Also considering
$$
\beta
= \frac{\frac{f(z_1)}{z_1}+\frac{f(z_2)}{z_2}+ \ldots +\frac{f(z_m)}{z_m}}{m}
$$
for $z_1,z_2, \ldots ,z_m \in \partial \mathbb{U}$,
we obtain

\

\begin{cor} \label{d4cor5.2} \quad
If $f(z)\in\mathcal{A}_{n}$ satisfies
$$
\left| \frac{zf'(z)}{f(z)}-1 \right|
< \frac{|1-\beta|n\rho}{1+|1-\beta|\rho}
\qquad(z\in\mathbb{U})
$$
for some complex $\beta=\dfrac{F(z_1)+F(z_2)+ \ldots +F(z_m)}{m} \in F(\mathbb{U})$
such that \\
$z_1,z_2, \ldots ,z_m \in \partial \mathbb{U}$ and $\beta\neq1$,
and for some real $\rho>1$
where $F(z)=\dfrac{f(z)}{z}$,
then
$$
\left| \frac{f(z)}{z}-1 \right|<\rho|1-\beta|
\qquad(z\in\mathbb{U}).
$$
\end{cor}

\


\begin{thebibliography}{}

\bibitem{m1ref1}
I. S. Jack,
{\it Functions starlike and convex of order $\alpha$},
J. London Math. Soc. {\bf 3}(1971), 469--474.

\bibitem{m1ref2}
S. S. Miller and P. T. Mocanu,
{\it Second-order differential inequalities in the complex plane},
J. Math. Anal. Appl. {\bf 65}(1978), 289--305.

\bibitem{m1ref0}
H. Shiraishi and S. Owa,
{\it Some sufficient problems for certain univalent functions},
Far East J. Math. Sci. {\bf 30}(2008), 147--155.

\end{thebibliography}
\end{document}